\documentclass[reqno]{amsart}
\usepackage{amscd,amsfonts,amssymb}
\usepackage{graphicx}
\usepackage{color}
\numberwithin{equation}{section}
\newtheorem{Theorem}{Theorem}[section]
\newtheorem{Lemma}{Lemma}[section]
\newtheorem{Corollary}{Corollary}[section]
\theoremstyle{definition}

\theoremstyle{remark}

\begin{document}

\begin{center}
{\bf Development of singularities for the compressible Euler
equations with external force in several dimensions}

{\it Olga Rozanova}

Department of Differential Equations,
Mathematics and Mechanics Faculty     \\
Moscow State University                 \\
Glavnoe zdanie GSP-2
Vorobiovy Gory
119992 Moscow
Russia \\
E-mail: {\it rozanova@mech.math.msu.su}
\vskip1cm
{\bf Abstract.}
\end{center}
{We consider solutions to the Euler equations in the whole space
from a certain class, which can be characterized, in particular,
by finiteness of mass, total energy and momentum. We prove that
for a large class of right-hand sides, including the viscous term,
such solutions, no matter how smooth initially, develop a
singularity within a finite time.
 We find a  sufficient condition
for the singularity formation,  "the best sufficient condition",
in the sense that one can explicitly construct a global in time
smooth solution for which this condition is not satisfied
"arbitrary little".

Also compactly supported perturbation of nontrivial constant state
is considered. We generalize the known theorem \cite {Sideris} on
initial data resulting in singularities. Finally, we investigate
the influence of frictional damping and rotation on the
singularity formation.}

\vskip2cm

\section{Introduction}

We are interested in the following system of balance laws
in differential form
$$ \rho (\partial_t
{\bf V}+({\bf V},{\bf \nabla})\,{\bf V})+{\bf\nabla} P=   {\bf
F}(t, x, \rho, {\bf V}, S, D^{|\alpha|}{\bf V}),\eqno(1)$$
$$ \partial_t \rho +
{\rm div}\, (\rho{\bf V})=0,\eqno(2)$$
$$ \partial_t S +({\bf
V},{\bf\nabla}S)=0,\eqno(3)$$ written for unknown functions $\rho,
{\bf V}=(V_1,...,V_n)$ and $S,$ density, velocity vector and
entropy, respectively. The functions depend on time $t$ and on
point $x=(x_1,...,x_n)\in {\mathbb R}^n.$ Here $P=P(t,x)$ is the
pressure, ${\bf F}=(F_1,...,F_n)$ is an external force, assumed to
be a smooth function of all its arguments, $|\alpha|\ge 0$ is a
multiindex, $\gamma>1$ is the adiabatic exponent.

We consider (1 -- 3) together with the state equation
$$P=\rho^\gamma e^S.\eqno(4)$$

For smooth solutions equations (2),(3) and (4) imply
$$ \partial_t P +({\bf
V},{\bf\nabla} P)+\gamma P\, div\,{\bf V}=0.\eqno(5) $$

Set an initial-value problem for (1), (2), (3), namely
$$\rho(0,x)=\rho_0(x), \,
{\bf V}(0,x)={\bf V}_0(x),\,S(0,x)=S_0(x).\eqno(6)$$ Sometime it
will be more convenient for us to consider the Cauchy problem for
(1),(2),(5), that is
$$\rho(0,x)=\rho_0(x), \, {\bf V}(0,x)={\bf
V}_0(x),\,P(0,x)=P_0(x)=\rho_0^\gamma(x) e^{S_0(x)}.\eqno(6')$$

It is well known that for at least ${\bf F}=0,$  solutions of
equations (1--3), no matter how smooth initially, can develop
singularities within a finite time.

In the one-dimensional (in space) case for the problem on a
singularity formation for solutions to (1--3)  the characteristics
method can be applied. In the isoentropic case, where the system
can be written in the Riemann invariants \cite{RozhdYanenko}), the
characteristics method gives a complete answer whether the
singularity (the gradient catastrophe) arises (it follows, for
example, from \cite{Borovikov}). In the non-isoentropic case also
there are some advances ( for example, \cite{Liu}), however the
results either have inexplicit character or concerns with small
perturbations of a constant state. The problem on the singularity
formation for the one-dimensional system of gas dynamic equations
can be investigated by means of other methods (\cite{Chemin,
Pohozhaev}). However, this methods give only sufficient conditions
of the gradient catastrophe, (generally speaking, with a large
margin).

In \cite{Sideris} for the 3D case sufficient conditions on initial
data perturbed from constant state with a positive density inside
a compact domain were found, such that the respective solution to
the Cauchy problem loses its smoothness in a finite time. The
results can be partly generalized to the case of arbitrary
dimension. The general sense of these conditions is that the speed
of the support propagation (that is, the  speed of sound in the
unperturbed domain) is small, compared with the velocity of the
gas inside  the initial perturbation. Provided the sufficient
conditions hold a breakdown  occurs  near the support boundary
\cite{Sideris}.

In \cite{Sideris} it is essential for the proof that the speed of
propagation of the perturbation is finite. Therefore the result
cannot be extended to the case of viscous compressible flow, with
the traditional viscosity description \cite{LL}, where the
perturbation spreads with infinite speed.

Problems where the initial data are compactly supported can be
treated separately. For this class of initial data it is
significantly easier to find initial conditions producing
singularities. The point is that if the solution is smooth, the
boundary of perturbation does not move, that is the support does
not expand. This fact allows to demonstrate that any smooth
compactly supported initial data result in a singularity
(\cite{Makino}). Moreover, it is true for the Navier-Stokes
equations as well (\cite{Xin}).

In Section 2 we consider initial data without restrictions on the
support, but having finite moment of mass and total energy. For
these solutions the mass is conserved. If we impose some
reasonable restrictions on the right-hand side of (1), then we
obtain additionally  conservation of angular momentum and
non-increasing of total energy. If the flow is considered in all
the space ${\mathbb R}^n,$ rather than just inside the bounded
volume of the liquid, then the density is forced to vanish rapidly
as $x\to\infty.$

We will show that in this case, too, it is possible to indicate
sufficient conditions to initial data, such that the solution
leaves a special class of functions. For some important right-hand
sides $\bf F$ it signifies that the solution loses its initial
smoothness within a finite time. The role of "restraining force"
preventing  decay of the gradient, rather than by the finite speed
of support propagation, is played by a value, that in the 3D case
can be interpreted as the initial vorticity of the flow.

This result is essentially multidimensional, in the sense that in
the 1D case the sufficient conditions  cannot be satisfied.

It is interesting that the result can be applied to the case of
viscous compressible flow, after imposing some restrictions on the
velocity vanishing at $|{\bf x}|\to \infty.$

Further we will show that the sufficient conditions that we find
are in some sense "best sufficient conditions" for a class of
right hand sides. That is, there exists an explicit globally
smooth in time  solution, for which the sufficient conditions are
not satisfied "arbitrarily little".

In Section 3 we improve the result of (\cite{Sideris}), and
generalize it, assuming the presence of a special exterior force,
that may have influence on the speed of the support propagation.

In Sections 4 and 5 we add to the right-hand side of (1) terms
describing damping and rotation and  find once more  sufficient
conditions for the finite time singularity formation.

\section{Solutions with finite moment of mass}
{\sc Definition.}{\it We will say that a solution $(\rho,{\bf
V},P)$ to system (1),(2),(5) belongs to the class $\mathfrak K$ if
it has the following properties:

(i)\quad the solution is classical for all $t\ge 0$;

(ii)\quad  $\rho|x|^2,\,P,\,\rho|{\bf V}|^2$ are of the class
$L_1({\mathbb R}^n);$

(iii)\quad $\displaystyle\int\limits_{{\mathbb R}^n}({\bf
F,x})\,dx\equiv 0,\quad $$\displaystyle\int\limits_{{\mathbb
R}^n}( F_i x_j-F_j x_i)\,dx\equiv 0,\,i\ne j,\,i,j=1,...,n,$ where
$\bf x$ is a radius-vector of point $x;$

(iv)\quad $\displaystyle\int\limits_{{\mathbb R}^n}({\bf
F,V})\,dx\le 0.$}

\medskip

We risk, of course, that for some choice of $\bf F$ the class
$\mathfrak K$ is empty  or consists only of trivial solution.
However, for $\bf F=0$ this class is not trivial and we
essentially seek the situation where the solutions can be treated
likely to this principal case.

For  ${\bf F}={\bf F(V)}, $ such that ${\bf F(0)}={\bf 0},$it is
known that if
$$(\rho_0^{\frac{\gamma-1}{2}}, {\bf V}_0, S_0)\in H^m({\mathbb R}^n)
,\quad m>\frac{n}{2}+1,$$ then locally in time (1--3) has a unique
solution
$$(\rho^{\frac{\gamma-1}{2}}, {\bf V}, S)\in
\cap_{i=0}^1 C^i([0,T],H^{m-i}({\mathbb R}^n).$$ This result
follows, for example, from
 \cite{Majda}, if the system is symmetrized by means of a new
 variable $P^{(\gamma-1)/2\gamma)}$\cite{MUK}. The classical result
 cannot be applied immediately, as on the solutions with finite
 moment of mass the density is not separated from zero and the
 system is not uniformely strictly hyperbolic.

Besides, in \cite{Chemin} it was proved if the initial data are
from the class $H_{ul}^m({\mathbb R}^n),$ then there exists a
unique solution from $\cap_{i=0}^1
C^i([0,T],H^{m-i}_{loc}({\mathbb R}^n).$ Here $H_{ul}^m$ is a
subset in $H^{m}_{loc}$ such that for all $\phi\in C_0^\infty,$ if
$\phi_x(y)=\phi(x-y)$, then
$\displaystyle\sup\limits_{x\in{{\mathbb R}^n}}\|\phi_x
u\|_{H^m({\mathbb R}^n)}<\infty.$

For the right hand side describing viscosity, there exist results
on a local in time existence of smooth solution as well, f.e.
\cite{Nash}, \cite{VolpertKhudiaev}. In \cite{Nash} the author
proved the existence of classical solution, having the H\"older
continuous second derivatives with respect the space variables and
the first ones with respect the time. In \cite{VolpertKhudiaev}
the system of equations of viscous compressible fluid is
considered as a particular case of composite systems of
differential equations. The consideration is proceeded in the
Sobolev spaces $H^l$ with a sufficiently large $l.$ The uniqueness
of the problem was proved earlier in \cite{Serrin}.

 For an
arbitrary forcing we have to assume the existence of a local in
time solution of class $\mathfrak K.$

Let us note that for the solutions of class $\mathfrak K$ we have
conservation of mass $m=\displaystyle\int\limits_{{\mathbb
R}^n}\rho\,dx.$ Moreover, all integrals
$$M_k=M_{ij}=\displaystyle\int\limits_{{\mathbb R}^n}(V_j x_i-V_i
x_j)\rho\,dx,\,i>j,\,k=1,...,K,\,K={\rm C}_n^2,$$ are conserved.
At $n=1$ there are no integrals in this series, at $n=2$ there is
only integral $M_1=\displaystyle\int\limits_{{\mathbb R}^n}({\bf
V, x_\bot})\rho\,dx,\,$ where ${\bf x}_\bot=(x_2,-x_1);$ at $n=3$
the integrals $M_1,M_2, M_3$ correspond to  components of the
angular momentum $\displaystyle\int\limits_{{\mathbb R}^n}({\bf
V\times x})\rho\,dx.$

 The total
energy, $E(t),$ is a sum of its kinetic and potential components,
that is
$$E(t)=E_k(t)+E_p(t):=\displaystyle\frac{1}{2}\int\limits_{{\mathbb
R}^n}\rho|{\bf V}|^2\,dx+\frac{1}{\gamma-1}\int\limits_{{\mathbb
R}^n}P\,dx.$$ We get from (1), (2), (5)   that $\displaystyle
E'(t)=\int\limits_{{\mathbb R}^n}({\bf F,V})\,dx.$ Thus, in virtue
of (iv), $E(t)$ is a non-increasing function for the solutions of
class $\mathfrak K.$

Let us introduce a functional
$$G_\phi(t)=\int\limits_{{\mathbb
R}^n}\rho(t,x)\phi(|{\bf x}|)\,dx,$$ considered for such functions
$\phi(|{\bf x}|)\in C^2[0,+\infty),$ that the integral converges
and
$$I_{4,\phi}(t)=\lim_{R\to \infty}\int\limits_{S(R)}\rho {\bf V}
\phi(|{\bf x}|) \,dS(R)=0,$$
 where $S(R)$ is the $(n-1)$ - dimensional sphere of radius $R$.

We denote by ${\bf M}=(M_1,...,M_K)$ and  ${\bf
\sigma}=(\sigma_1,...,\sigma_K)$  vectors with  components
$M_k,\,$ and  $\sigma_k=V_i x_j-V_j x_i,
\,i>j,\,i,j=1,...,n,\,k=1,...,K,\,K={\rm C}_n^2,$ respectively.

\begin{Lemma} Provided all given integrals converge,
for  solutions to (1), (2), (5) of the class $\mathfrak K$
following equalities take place:
$$G'_\phi(t)=\int\limits_{{\mathbb R}^n}\frac{\phi'(|{\bf x}|)}{|{\bf x}|}
({\bf V},{\bf x})\rho\,dx,$$

$$G''_\phi(t)=I_{1,\phi}(t)+I_{2,\phi}(t)+I_{3,\phi}(t)+I_{4,\phi}(t),$$
where
$$I_{1,\phi}(t)=\int\limits_{{\mathbb R}^n}\frac{\phi''(|{\bf x}|)}{|{\bf
x}|^2} |({\bf V},{\bf x})|^2\rho\,dx,$$
$$I_{2,\phi}(t)=\int\limits_{{\mathbb R}^n}\frac{\phi'(|{\bf x}|)}{|{\bf
x}|^3} |{\bf \sigma}|^2\rho\,dx,$$
$$I_{3,\phi}(t)=\int\limits_{{\mathbb R}^n}(\phi''(|{\bf x}|)+
(n-1)\frac{\phi'(|{\bf x}|)}{|{\bf x}|} ) P\,dx,$$
$$I_{4,\phi}(t)=-\lim_{R\to \infty}\int\limits_{S(R)}\phi'(|{\bf x}|)
 P\,dS(R).$$
\end{Lemma}

The proof is a direct calculation and an application of the
general Stokes formula. For example, we get, using (2), that
$$G'_\phi(t)=\int\limits_{{\mathbb R}^n}\rho'_t(t,x)\phi(|{\bf
x}|)\,dx=-\int\limits_{{\mathbb R}^n}{\rm div}(\rho{\bf
V})\phi(|{\bf x}|)\,dx=$$$$=\int\limits_{{\mathbb
R}^n}(\nabla\phi(|{\bf x}|), {\bf V})\rho\,dx-\lim_{R\to
\infty}\int\limits_{S(R)}
 \rho{\bf V}\phi(|{\bf x}|)\,dS(R)=$$
$$=\int\limits_{{\mathbb R}^n}\frac{\phi'(|{\bf
x}|)}{|{\bf x}|} ({\bf V},{\bf x})\rho\,dx.$$ $\Box$

{\sc Remark 2.1} If the increase of $\phi(|{\bf x}|)$ as $|{\bf
x}|\to \infty$ is no more then $const\cdot |{\bf x}|^2,$ then the
condition
$$\lim_{R\to
\infty}\int\limits_{S(R)}
 \rho{\bf V}\phi(|{\bf x}|)\,dS(R)=0$$
 follows from  (ii) without additional assumptions on the
 behavior of velocity at
$|{\bf x}|\to \infty.$
\medskip

In the particular case $\phi(|{\bf x}|)=\frac{|{\bf x}|^2}{2}$ we
denote $G_\phi(t),I_{i,\phi}(t),\,i=1,...,4,$ by $G(t),I_{i}(t),$
respectively, the derivative $G'(t)$ we denote $F(t).$

\begin{Corollary}
For  solutions to (1), (2), (5) of the class $\mathfrak K$
$$F(t)=G'(t)=\int\limits_{{\mathbb R}^n}
({\bf V},{\bf x})\rho\,dx,$$
$$I_{1}(t)=\int\limits_{{\mathbb R}^n}\frac{|({\bf V},{\bf x})|^2
}{|{\bf x}|^2} \rho\,dx,$$
$$I_{2}(t)=\int\limits_{{\mathbb R}^n}\frac{|{\bf \sigma}|^2}
{|{\bf x}|^2} \rho\,dx,$$
$$I_{3}(t)=n\int\limits_{{\mathbb R}^n}P\,dx=n(\gamma-1)E_p(t),$$
$$I_{4}(t)=0.$$
Moreover,
$$I_{1}(t)+I_{2}(t)=2E_k(t).$$
\end{Corollary}

The proof of Corollary 2.1 is an immediate substitution of a
particular form of $\phi(|{\bf x}|);$ $\,I_{4}(t)=0$ due to
sufficiently rapid vanishing of $P$ as $|{\bf x}|\to\infty,$
forced by condition
 (ii). $\Box$

 The following Lemma gives
some useful estimates for the functionals introduced above.

\begin{Lemma}
For solutions to (1), (2), (5) of the class $\mathfrak K$
inequalities
$$(G'(t))^2=F^2(t)\le 4G(t)E_k(t)\le 4G(t)E(0),\eqno(7)$$
$$(G'(t))^2=F^2(t)\le 2G(t)I_{1}(t),\eqno(8)$$
$$|{\bf M}|^2\le 4G(t)E_k(t)\le 4G(t)E(0),\eqno(9)$$
$$|{\bf M}|^2\le 2G(t)I_{2}(t),\eqno(10)$$
$$G(t)\le (\sqrt{E(0)}t+\sqrt{G(0)})^2\eqno(11)$$
holds.
\end{Lemma}

Proof. The first four inequalities are  corollaries of the
H\"older inequality. Inequality (11) follows from (7) after
integration.$\Box$

\medskip

Let us point out that in the case ${\bf F=0}$ the last parts in
(7) and (9) is not a very strong roughening, as according to
\cite{Chemin}, for smooth solutions $E_k(t)\to E,\, t\to \infty$.

 The H\"older inequality gives us also a lower
estimate for the kinetic energy, namely,
$$E_k(t)\ge\frac{F^2(t)}{4G(t)}.\eqno(12)$$

As for a lower estimation of the potential energy, there exists
the following result.
\medskip
\begin{Lemma}  \cite{Chemin} For
solutions to (1), (2), (5) satisfying (ii)
$$E_p(t)\ge \frac{C}{G^{\frac{(\gamma-1)n}{2}}(t)},\eqno(13)$$
with a positive constant $C,$  depending on initial data, $\gamma$
and $n.$

If we denote $\displaystyle S_0=\inf\limits_{x\in {\mathbb
R}^n}S(0,x),$ then
$$C=\frac{e^{S_0}}{\gamma-1}(m C_{\gamma,
n}^{-1})^\frac{\gamma(n+2)-n}{2},$$ with $$C_{\gamma,
n}=\left(\frac{2\gamma}{n(\gamma-1)}\right)^
{\frac{n(\gamma-1)}{(n+2)\gamma-n}} +
\left(\frac{2\gamma}{n(\gamma-1)}\right)^
{\frac{-2\gamma}{(n+2)\gamma-n}}.
$$
\end{Lemma}

\medskip\medskip

\begin{Theorem}
There are initial data (6') satisfying (ii)  such that
 solution to (1), (2), (5) from the class $\mathfrak K$
 exists during a
finite time. Namely, it occurs if
$$  F(0)\ge L_2 \cot
\frac{L_1}{2\sqrt{E(0)G(0)}},\eqno(14)$$ with constants $L_1$ and
$L_2$ depending on initial data , the adiabatic exponent $\gamma$
and the dimension of space only.

 If  $\gamma\in(1,1+\frac{2}{n}],$ then
$L_1=L_2=L:=\left(2n(\gamma-1)C(G(0))^{1-(\gamma-1)n/2}+|{\bf
M}|^2 \right)^{1/2}.$

If $\gamma>1+\frac{2}{n},$ then $L_1=L/((\gamma-1)n-1),\,L_2=L.$

The time of existence for this solution can be estimated above by
the constant
$$T_*=\sqrt{\frac{G(0)}{E(0)}}\frac{
\frac{2\sqrt{E(0)G(0)}}{L_2}
\left(\frac{\pi}{2}-\arctan\frac{F(0)}{L_1}\right)}
{1-2\frac{\sqrt{E(0)G(0)}}{L_2}\left(\frac{\pi}{2}-\arctan\frac{F(0)}
{L_1}\right)}.
$$
\end{Theorem}

\medskip
Proof of Theorem 2.1. Let $\gamma\in(1,1+\frac{2}{n}],$ it follows
$2-(\gamma-1)n\ge 0.$ Therefore, from (8), (10), (13) we have
$$F'(t)\ge
\frac{F^2(t)}{2G(t)}+\frac{|{\bf
M}|^2}{2G(t)}+\frac{n(\gamma-1)C}{(G(t))^{\frac{(\gamma-1)n}{2}}}
\ge
$$
$$\ge
\frac{F^2(t)+|{\bf
M}|^2+2n(\gamma-1)C(\sqrt{E(0)}t+\sqrt{G(0)})^{2-(\gamma-1)n}}
{2(\sqrt{E(0)}t+\sqrt{G(0)})^{2}}\ge$$$$\ge
\frac{F^2(t)+L^2}{2(\sqrt{E(0)}t+\sqrt{G(0)})^{2}}.\eqno(15)$$
After integration we have
$$\arctan \frac{F(t)}{L}\ge \arctan \frac{F(0)}{L}
+\frac{L}{2\sqrt{E(0)}}\left(\frac{1}{\sqrt{G(0)}}-
\frac{1}{\sqrt{E(0)}t+\sqrt{G(0)}}\right).\eqno(16)$$ The left
hand side of (16) does not exceed $\frac{\pi}{2},$ therefore (16)
cannot be true for all $t>0$ if
$$\arctan \frac{F(0)}{L}
+\frac{L}{2\sqrt{E(0)G(0)}}> \frac{\pi}{2}.\eqno(17)$$

Let us note that
$$\frac{L}{2\sqrt{E(0)G(0)}}\le\frac{\sqrt{2n(\gamma-1)E_p(t)G(0)+4E_k(t)
G(0)}}{2\sqrt{(E_p(0)+E_k(0))G(0)}}\le 1<\frac{\pi}{2},$$
therefore after trigonometric transformations of (17) one can get
(14).

Further, let
$\gamma >1+\frac{2}{n},$ therefore $(\gamma-1)\frac{n}{2}-1\ge 0.$
Analogously to the previous case
we have
$$F'(t)\ge
\frac{F^2(t)+|{\bf
M}|^2}{2G(t)}+\frac{n(\gamma-1)C}{(G(t))^{\frac{(\gamma-1)n}{2}}}
\ge
$$
$$\ge
\frac{(F^2(t)+|{\bf
M}|^2)(\sqrt{E(0)}t+\sqrt{G(0)})^{-1+(\gamma-1)n}+ n(\gamma-1)C}
{2(\sqrt{E(0)}t+\sqrt{G(0)})^{(\gamma-1)n}}\ge
$$$$\ge\frac{(G(0))^{\frac{(\gamma-1)n}{2}-1}(F^2(t)+|{\bf
M}|^2)+2n(\gamma-1)C}
{2(\sqrt{E(0)}t+\sqrt{G(0)})^{(\gamma-1)n}}=$$
$$=
\frac{F^2(t)+L^2} {2(G(0))^{\frac{-(\gamma-1)n}{2}+1}
(\sqrt{E(0)}t+\sqrt{G(0)})^{(\gamma-1)n}}. \eqno(18)$$ After
integration we get
$$\arctan \frac{F(t)}{L}\ge \arctan \frac{F(0)}{L}
+
$$
$$\frac{L}{2((\gamma-1)n-1)\sqrt{E(0)}(G(0))^{\frac{-(\gamma-1)n}{2}+1}}
\left(\frac{1}{(G(0))^{\frac{(\gamma-1)n-1}{2}}}-
\frac{1}{(\sqrt{E(0)}t+\sqrt{G(0)})^{(\gamma-1)n-1}}\right).\eqno(19)$$
The condition (19) cannot be true for all $t>0$ if
$$\arctan \frac{F(0)}{L}
+\frac{L}{2((\gamma-1)n-1)(\sqrt{E(0)G(0)}}\ge \frac{\pi}{2}.$$ It
implies (14) analogously to the previous case. Theorem 2.1 is
proved.$\Box$

\medskip
\medskip

{\sc Remark 2.2.} It seems that we can obtain an analogous
nonexistence result from (8) and Lemma 2.1 using only the
nonnegativity of the integrals $I_{2}(t)$ and $I_{3}(t)$. However,
it is not true. Indeed, here we have
$$F'(t)\ge
\frac{F^2(t)}{2G(t)}\ge \frac{F^2(t)}
{2(\sqrt{E(0)}t+\sqrt{G(0)})^{2}}.\eqno(20)$$ So we obtain
$$-\frac{1}{F(t)}+\frac{1}{F(0)}\ge
\frac{1}{2\sqrt{E(0)}}\left(\frac{1}{\sqrt{G(0)}}-
\frac{1}{\sqrt{E(0)}t+\sqrt{G(0)}}\right).$$ As $F'(t)>0,$ then
$F(t)$ remains positive for $F(0)>0.$ Therefore we conclude that
(20) cannot hold for all $t$ if
$$F(0)> 2\sqrt{E(0)G(0)}.\eqno(21)$$
However, (21) contradicts the inequality (7), therefore we cannot
choose the initial data with such properties.

It is interesting that if ${\bf M=0},$ then one cannot find
initial data satisfying (14), either.

To show this, let us consider, for example, the case $\gamma \in
(1,1+\frac{2}{n}].$ Let us find a necessary condition for the
implementation of (14). As follows from (7) and (14)
$$\frac{2\sqrt{E_k(0)G(0)}}{L}\ge
\cot\frac{L} {2\sqrt{(E_k(0)+E_p(0))G(0)}}.\eqno(22)
$$
We denote $\displaystyle
z=\frac{2\sqrt{E_k(0)G(0)}}{L},\,z_1=\frac{2\sqrt{E_p(0)G(0)}}{L}.$
Further,  we introduce a function
$$f(z):=\arctan \frac{1}{z}-\frac{1}{\sqrt{z^2+z_1^2}},
\quad z\in [0,\infty). \eqno(23)$$ Since  (22) signifies
$$\frac{1}{z}<\tan \frac{1}{\sqrt{z^2+z_1^2}},$$
 then for implementation of
condition (22) we have to find a point $z_*$ such that $f(z_*)\le
0.$ However due to Lemma 2.3 for ${\bf M=0}$ we have
$$z_1=\frac{2\sqrt{E_p(0)G(0)}}{\sqrt{2n(\gamma-1)CG^{1-(\gamma-1)n/2}(0)
}} \ge \frac{2\sqrt{E_p(0)G(0)}}{\sqrt{2n(\gamma-1)E_p(0)G(0)}}\ge
1.$$ Therefore $f(0)=\frac{\pi}{2}-\frac{1}{z_1}>0,$ and one can
show that $f(z)$ will be positive for all $z>0.$

On the other side, if $E_p(0)=0$ (it take place in the so called
"pressureless" gas dynamic, when $P\equiv 0$ \cite{Zeldovich}),
then (14) can be satisfied  also for ${\bf M=0}.$

From the consideration above we can conclude that the condition
(14) cannot hold in one space dimension if $E_p(0)\ne 0$, where we
cannot obtain the additional positive lower bound for the integral
$I_{2}(t).$

 As follows from (22), the necessary
condition of implementation of (14) is the negativity of $f(z)$ at
some points. For ${\bf M}\ne 0$, then it will be, for example, if
$z_1\le \frac{\pi}{2},$ that is $L\ge {\pi}{\sqrt{E_p(0)G(0)}}.$
The last inequality surely holds if $|{\bf
M}|>\frac{\sqrt{E_p(0)G(0)}}{\pi}.$ It implies
$$\frac {E_k(0)}{E_p(0)}\ge \frac{\pi^2}{4},$$ that is initially
the part of kinetic energy must exceed the potential one.

Now we will show that together with a large value of $|{\bf M}|,$
the condition (14) requires a large initial divergency of the
flow. We denote now $Z:=\frac{L}{2\sqrt{EG(0)}}$ and point out
that $Z\le 1.$ Condition (14) can be re-written as
$$\lambda Z\le \tan Z,$$
with $\lambda=\frac{2\sqrt{EG(0)}}{F(0)}$, therefore $\lambda<\tan
1,$ or $F(0)\ge \frac{2\sqrt{EG(0)}}{\tan 1}.$ It follows from the
last inequality that $\frac{E_k(0)}{E}\ge \cot 1.$

\medskip

{\sc Remark 2.3.} As follows from \cite{Chemin1}, the breakdown of
smoothness in the compressible non-viscous flow in 3D comes from
the accumulation of vorticity, divergency or compression.

\medskip

{\sc Remark 2.4.} As follows from the proof of Theorem 2.1, the
singularity appearance is a result of an unlimited growth of
$F(t)$. For solutions from the class $\mathfrak K$ the Green's
formula shows that
$$F(t)=-\frac{1}{2}\int\limits_{{\mathbb R}^n}|{\bf x}|^2\nabla(\rho{\bf
V})d{\bf x},$$ therefore the predicted appearence of a singularity
can be associated with domains of "large negative divergency" or,
as meteorologists say, of  "large convergency."

In its turn, in the physical space 3D $$|{\bf
M}|=\left|\frac{1}{2}\int\limits_{{\mathbb R}^3}|{\bf x}|^2 {\rm
rot}(\rho{\bf V})d{\bf x}\right|,$$ that is a large value of
$|{\bf M}|$ corresponds to a large initial vorticity.

\medskip

\subsection{On  exactness of integral condition (14)}

Since there is no one-to-one correspondence between solutions to
the system (1), (2), (5) and the integral functionals considered,
we cannot expect that the condition (14) of Theorem 2.1 are
sufficient and necessary conditions for the singularity appearance
in the class of solutions with a finite total energy and a finite
moment of mass satisfying (iii) and (iv).

Nevertheless, further we will see that condition (14) of Theorem
2.1 is "exact sufficient condition" for singularity appearance, at
least for special right-hand sides. In other words, if it is not
satisfied for certain initial data "arbitrary little", then the
corresponding solution to the Cauchy problem may be globally
smooth in time. Namely, the following Theorem holds:

\begin{Theorem} Let $\bf F =0 $ for any velocity field of form
${\bf V}=\alpha(t){\bf x}.$ Then for an arbitrary small
$\varepsilon>0$ there exists a globally in time classical solution
to system (1), (2), (5) from the class $\mathfrak K$, such that at
$t=0$ condition
$$ F(0)>L_1 \cot
\frac{L_2}{2\sqrt{E(0)G(0)}}-\varepsilon\eqno(24)$$ holds.
\end{Theorem}

\medskip
To prove  Theorem 2.2 we give an example of solutions satisfying
the properties indicated in the Theorem statement. It is known
that there exists a class of globally smooth solutions with linear
profile of velocity  (for one-dimensional case one can find its
description in \cite{Sedov},\cite{Ovs}, where the Lagrangean
variables are used; another approach to constructing and
generalization to the case of several space dimensions there are
in \cite{RozPas}, \cite{RozJMS}, \cite{RozNova}). It will be
sufficient for us to consider the simplest form of such fields of
velocity, namely,
$$
{\bf V}=\alpha(t){\bf x},\eqno(25)
$$
with  a function $\alpha(t),$ taking part of the solution  to
system of ODE

$$G_1'(t)=-2\alpha(t) G_1(t),\quad \alpha'(t)=-\alpha^2(t)+
(\gamma-1)K G_1^{\frac{(\gamma-1)n}{2}+1}(t),\eqno(26)
$$

Here $K={ E}_p(0) G^{\frac{(\gamma-1)n}{2}}(0),$$\,G_1(t)=
1/G(t).$

The components of density and pressure can be found from the
linear with respect to them equations (2) and (5) as
$$
\rho(t,|{\bf x}|, \phi)=\exp(-2\int\limits_0^t \alpha(\tau) d\tau)
\rho_0(|{\bf x}|\exp(-\int\limits_0^t \alpha(\tau) d\tau)), $$
$$
p(t,|{\bf x}|, \phi)=\exp(-2\gamma\int\limits_0^t \alpha(\tau)
d\tau) p_0(|{\bf x}|\exp(-\int\limits_0^t \alpha(\tau) d\tau)), $$
with compatible initial data $\rho_0(x),P_0(x).$ The compatibility
signifies here that the condition
$${\bf
\nabla} p_0(x)= - (\gamma-1)G_1(0)E_p(0)\rho_0(x) {\bf x}\eqno
(27)$$ holds.

For example, one can choose
$$p_0=\frac{1}{(1+|{\bf
r}|^2)^a},\,a=const>\frac{n}{2}, $$
$$ \rho_0=
\frac{2a}{(\gamma-1)G_1(0)E_p(0)}\frac{1}{(1+|{\bf
r}|^2)^{a+1}}.$$

Let us note that system (26) takes place for all solutions with
the velocity profile (25), however, compatibility condition (27),
generally speaking, can do not hold, therefore we have to require
the special form of ${\bf F}.$

For the solutions considered $F(t)=2\alpha(t)G(t),\quad {\bf
M}=0,$
 the
kinetic energy $E_k(t)=\alpha^2(t)G(0),$ the potential energy
$E_p(t)=\displaystyle\frac{K}{(G(t))^{(\gamma-1)n/2}}.$

Let us consider, for example, the case $\gamma\le 1+\frac{2}{n}.$

 Thus, (24) takes the
form
$$\frac{2\alpha(0)
G(0)}{\sqrt{2n(\gamma-1)CG(0)^{1-(\gamma-1)n/2}}}\ge
\cot\frac{\sqrt{2n(\gamma-1)CG(0)^{1-(\gamma-1)n/2}}}{2\sqrt{G(0)}
(\alpha^2(0)G(0)+KG^{-(\gamma-1)n/2})}-\varepsilon$$ or

$$\frac{2\alpha(0)
(G(0))^{\frac{1}{2}+\frac{(\gamma-1)n}{4}}}{\sqrt{2n(\gamma-1)C}}\ge
\cot\frac{\sqrt{2n(\gamma-1)C}}{2\sqrt{
\alpha^2(0)(G(0))^{1+(\gamma-1)n/2}+K}}-\varepsilon.\eqno(28)$$

Let us fix $\rho_0(x)$ and $P_0(x).$ It signify that $G(0)$ and
$C$ are fixed. We are going to show that we can choose $\alpha(0)$
such that for anyhow small positive $\varepsilon$ inequality (28)
will be satisfied.

We denote
$$z(\alpha(0))=\frac{1}{\alpha(0)}\frac{\sqrt{2n(\gamma-1)C}}{
2(G(0))^{\frac{1}{2}+\frac{(\gamma-1)n}{4}}}$$ and
$$\lambda(\alpha(0))=\frac{\alpha(0)(G(0))^{\frac{1}{2}+
\frac{(\gamma-1)n}{4}}}{\sqrt{
\alpha^2(0)(G(0))^{1+(\gamma-1)n/2}+K}}.$$ We note that
$z(\alpha(0))\to 0$ and $\lambda(\alpha(0))\to 1$ as
$\alpha(0)\to\infty,$ moreover, $\lambda(\alpha(0))<1$ for any
finite $\alpha(0).$

Thus, (28) can be re-written as follows:
$$\frac{1}{z[(\alpha(0)]}\ge\cot z[(\alpha(0)]+\left(\cot
[\lambda(\alpha(0))z(\alpha(0))]-\cot
z[\alpha(0)]\right)-\varepsilon.\eqno(29)$$ We point out that if
$z(\alpha(0))\in (0,\pi),$ then $\cot
[\lambda(\alpha(0))z(\alpha(0))]-\cot [z(\alpha(0))]>0,$ however,
for any $\varepsilon>0$ we can choose $\alpha_0>0$ such that for
any $\alpha(0)>\alpha_0$ the difference $\left(\cot
[\lambda(\alpha(0))z(\alpha(0))]-\cot
z[\alpha(0)]\right)-\varepsilon<-\varepsilon_1$ for some
$\varepsilon_1>0.$

Since $\displaystyle\frac{1}{z}\ge\cot z,\,z\in (0,\pi),$ then
$\displaystyle\frac{1}{z}\ge\cot z-\varepsilon_1.$  For
$\alpha(0)>\alpha_0$ it implies  (29) and, consequently, (24).

The case $\gamma>1+\frac{2}{n}$ can be treated analogously.

Thus, the proof of Theorem 2.2 is over.$\Box$

\medskip
{\sc Remark 2.5.} Besides the trivial case ${\bf F=0}$, the first
condition of Theorem 2.2 is satisfied for ${\bf
F=F}(D^{|\alpha|}{\bf V}),\,|\alpha|\ge 2.$

\medskip
{\sc Remark 2.6.} However, we cannot assert that if condition (14)
does not hold for certain initial data, then the solution to the
corresponding Cauchy problem is necessarily globally smooth in
time. For example, for ${\bf F=0},$ let us consider initial data
with zero velocity and compactly supported density (and pressure).
These initial data always result in a singularity
(f.e.\cite{Makino},\cite{RozMag}). However, $F(0)=0,$ and as the
right hand side in (14) is positive, the condition (14) of Theorem
2.1 is not satisfied. On the other hand, it seems that in this
situation there exists a moment $t_1$ such that if it is chosen as
the initial one, then (14) will be already satisfied. In others
words, the hypothesis is that (14) detects singularities arising
from accumulation of negative divergency, which are sufficiently
close in time.

\medskip

{\sc Remark 2.7.} One should pay attention to the following fact:
smooth initial data $\rho_0(x),P_0(0),{\bf V}(0),$ having compact
support and satisfying the compatibility condition are not good
for application of the theorem on a local in time existence and
uniqueness of the Cauchy problem for the symmetric hyperbolic
systems (\cite{Majda}). The matter is that at the point where the
density vanishes smoothly, for a compatible initial data the
entropy becomes infinite,  so we cannot apply the cited theorem,
which require the smoothness of initial data for {\it symmetrized}
system, where the variables are entropy, velocity and
$P^{(\gamma-1)/2\gamma}$ (see \cite{Makino}). Indeed, we get non-
uniqueness for the compatible initial data of density and pressure
as follows. According to  procedure described in the Theorem 2.2
proving we can construct a global in time solution with the
velocity field of form $\alpha(t){\bf x}.$ Let us choose a moment
$t_0$ such that $\alpha(t_0)=0.$ Then the initial velocity ${\bf
V}_0(x)=0.$ At the same time it is known that the solution with
smooth density having a compact support cannot be globally smooth.
Note that if the density and pressure are only continuous at the
points of vanishing, then they can be compatible. On can construct
this solution; its support spreads.
\medskip

\subsection{Application to the compressible Navier-Stokes
system}

For the Navier-Stokes system, describing the behavior of
compressible viscous fluid, the right hand side of (1) is the
following:
$${\bf F}={\rm div} T, \quad T_{ij}=\mu(\partial_i V_j+
\partial_j V_i) +\lambda {\rm div}V
\delta_{ij},\,i,j=1,...,n,\eqno(30)$$ where $T$ is the stress
tensor, $\mu\ge 0,$ and $\lambda$ are constants
$(\lambda+\frac{2}{n}\mu\ge 0),\,\delta_{ij}$ is the Kronekker
symbol.

In \cite{Xin} it was demonstrated that if
$\mu>0,\,\lambda+\frac{2}{n}\mu> 0,$  then there exists no
solution with compactly supported density to the Cauchy problem
(1), (2), (5), (6') with the right-hand side of form (30) from
$C^1([0,\infty), H^m({\mathbb R}^n)),\quad
m>2+\left[\frac{n}{2}\right]$, such that initial data (6') are in
$ H^m({\mathbb R}^n).$

As a corollary of Theorem 2.1 we obtain that if the density and
velocity  vanish at infinity sufficiently quickly, there are
initial conditions such that the solution to the Cauchy problem
exists only within a finite interval of time.

Let $S_R$ be a sphere of radius $R,$ the unit outer normal to
$S_R$ and the element of its surface we denote by ${\bf
N}(N_1,...,N_n)$ and $d S_R, $ respectively.
\medskip

{\sc Definition} {We will say that a solution to (1),(2), (5) with
the right-hand side of form (30) is of class ${\mathfrak K}_1$ if
it satisfies to conditions (i),(ii) and
$$(v)\quad\lim\limits_{R\to\infty}\int\limits_{S_R}
\sum\limits_{i,j=1}^n (T_{ij}x_i
-(2\mu+n\lambda)V_i\delta_{ij})N_j d S_R=0,$$
$$(vi)
\lim\limits_{R\to\infty}\int\limits_{S_R} \sum\limits_{i,j=1}^n
T_{ij}V_i N_j d S_R=0.$$}

\medskip

{\sc Remark 2.8.} Conditions (v),(vi) are satisfied, for example,
if the velocity vector decays as $|{\bf x}|\to \infty$ uniformly
in $t$ so that $|{\bf V}|=o\left(\frac{1}{|x|^{n-1}}\right)$ and
$|D{\bf V}|=o\left(\frac{1}{|x|^n}\right).$

\begin{Theorem}
For initial data (6') satisfying  condition(14) of Theorem 2.1,
the solution to (1), (2), (5) from the  class ${\mathfrak K}_1$
exists only within a finite time.
\end{Theorem}

To prove the Theorem 2.3 it suffices to note that in this
situation ${\mathfrak K}_1\subset \mathfrak K,$ because the rate
of the velocity decay implies
$$\displaystyle\int\limits_{{\mathbb R}^n}({\bf F,x})\,dx
=\int\limits_{{\mathbb
R}^n}\sum\limits_{i,j=1}^n\partial_jT_{ij}x_i dx =
\lim\limits_{R\to\infty}\int\limits_{S_R} \sum\limits_{i,j=1}^n
T_{ij}x_i N_j d S_R -(2\mu+n\lambda)\int\limits_{{\mathbb R}^n}div
{\bf V}\,dx= $$ $$=\lim\limits_{R\to\infty}\int\limits_{S_R}
\sum\limits_{i,j=1}^n (T_{ij}x_i
-(2\mu+n\lambda)V_i\delta_{ij})N_j d S_R=0,$$

$$\displaystyle\int\limits_{{\mathbb R}^n}(F_i x_j-F_j x_i)\,dx
= \lim\limits_{R\to\infty}\int\limits_{S_R}
\sum\limits_{i,j,k=1}^n (T_{kj}x_i- T_{ki}x_j)N_k d S_R= 0,$$

$$\displaystyle\int\limits_{{\mathbb R}^n}({\bf F,V})\,dx
=\int\limits_{{\mathbb
R}^n}\sum\limits_{i,j=1}^n\partial_jT_{ij}V_i dx = $$
$$
=\lim\limits_{R\to\infty}\int\limits_{S_R} \sum\limits_{i,j=1}^n
T_{ij}V_i N_j d S_R- \int\limits_{{\mathbb
R}^n}\sum\limits_{i,j=1}^nT_{ij}\partial_i V_j dx \le 0.$$
$\Box$
\medskip

We point out that classes $\mathfrak K$ and ${\mathfrak K}_1$ do
not coincide. For example, the solutions with linear profile of
velocity belong to $\mathfrak K$, however, condition (vi) is not
satisfied here and therefore these solutions are not of class
${\mathfrak K}_1.$

\medskip

{\sc Remark 2.9.} The fact that a solution with a finite moment of
mass and a finite total energy leaves the class $\mathfrak K$ both
in the forcing free Euler system and the Navier-Stokes system
means that a singularity appears (provided the decay of velocity
at infinity in the Navier-Stokes case is sufficiently rapid). But
the nature of the singularity appearing under the same initial
conditions for these systems is different. In both situations, it
signifies that the integration by parts becomes prohibited.
However, in the first case this integration is still allowed if
the solution is only continuous along certain piece-wisely smooth
curves. For Navier-Stokes system we need to require  $ C^1 $ --
smoothness with respect to the space variables along these curves.
Thus, for the hyperbolic systems the  singularity predicted is
either a strong discontinuity or some week discontinuity on a
complicated set. For the Navier-Stokes system even week
discontinuity along smooth curves means that a singularity
appears.

On the other hand, in the case of the Navier-Stokes system "the
singularity appearance" may signify that the solution does not
belong anymore to the class ${\mathfrak K}_1.$ For example, the
derivatives of velocity in this possible global-in-time smooth
solution satisfying (14) cannot decay at infinity sufficiently
quickly.
\medskip

{\sc Remark 2.10.} Condition (14) of Theorem 2.1 is "the best
sufficient condition" for the leaving the class ${\mathfrak K},$
in the viscous case, too. Indeed, the viscous term "does not feel"
the velocity with a linear profile, therefore the first condition
of Theorem 2.2 is satisfied, and we can apply Theorem 2.2 in this
situation. It is interesting that this solution belongs to
${\mathfrak K}\backslash{\mathfrak K}_1.$

\section{Breakdown of compactly supported smooth
perturbation of a constant state}

Results of \cite {Sideris} for the perturbation of the constant
state $(\bar\rho, {\bf 0}, \bar S)$ having compact support $B(t)$
can be extended to the case of right-hand sides ${\bf F}$ with the
properties, analogous to (iii) and (iv), if we suppose the finite
speed of the perturbation propagation. This characteristic
property of hyperbolic systems, generally speaking, does not hold
for Navier-Stokes equations (at least, for $\bar\rho\ne 0,$ see in
this context \cite{Xin}). Indeed, to obtain  results analogous to
\cite {Sideris}, it is sufficient to impose certain condition to
the function $R(t),$ where $R(t)$ is the minimal radius of ball
containing the support of perturbation. This condition has the
form:
$$R(t)< C (1+t)^\alpha,\, \alpha \in {\mathbb
R},\, C\in {\mathbb R}_+.\eqno(31)$$ If ${\bf F}={\bf 0},$ then
$\alpha=1.$

Let us denote, following to \cite {Sideris},
$$m(t)=\int\limits_{{\mathbb R}^n}[\rho(t,x)-\bar\rho] dx=
\int\limits_{B(t)}[\rho(t,x)-\bar\rho] dx,$$
$$
\eta(t)=\int\limits_{{\mathbb
R}^n}[\rho(t,x)\exp\left(\frac{S(t,x)}{\gamma}\right)-
\bar\rho\exp\left(\frac{\bar S}{\gamma}\right)]
dx=$$$$=\int\limits_{B(t)}[\rho(t,x)\exp\left(\frac{S(t,x)}{\gamma}\right)-
\bar\rho\exp\left(\frac{\bar S}{\gamma}\right)] dx,$$
$$\tilde G(t)=\frac{1}{2}\int\limits_{B(t)}|{\bf x}|^2\rho dx,
\quad \tilde F(t)=\int\limits_{B(t)}({\bf x,V})\rho dx.
$$
If we recall the denotation of section 2, it occurs that
$F(t)=\tilde F(t)$ (see Corollary 2.1).

We obtain the following generalization of Theorem 1 from \cite
{Sideris}:

\begin{Theorem}
Let us suppose that $(\rho, V, P)$ is a classical solution to
(1,2,5) such that  properties

\begin{itemize}

\item $\displaystyle\int\limits_{B(t)}({\bf F,x})\,dx\equiv 0;$
\item $\displaystyle\int\limits_{B(t)}( F_i x_j-F_j x_i)\,dx\equiv
0,\,i\ne j,\,i,j=1,...,n,$ \noindent where $\bf x$ is a
radius-vector of the point $\,x;$ \item
$\displaystyle\int\limits_{B(t)}({\bf F,V})\,dx\le 0$
\end{itemize}
hold.

We assume that the support of perturbation propagates according to
condition (31)and $\eta(0)\ge 0.$  Suppose also that any of
following conditions takes place:

a) \quad  $\alpha\le\frac{1}{2+n},\quad|{\bf M}|=0, \quad F(0)>0;$

b)\quad  $\alpha\le\frac{1}{2+n},\quad|{\bf M}|\ne 0;$

c)\quad $\alpha>\frac{1}{2+n},\,$  $\,F(0)> (\alpha(2+n)-1)A,$

d)\quad $\alpha>\frac{1}{2+n},\,$  $\displaystyle\,F(0)> |{\bf
M}|\cot \frac{|{\bf M}|}{(\alpha(2+n)-1)A}, \,|{\bf M}|\ge \pi
(\alpha(2+n)-1) A,$

e)\quad $\alpha>\frac{1}{2+n},\,$  $|{\bf M}|\ge \pi
(\alpha(2+n)-1) A,$

where the constant $A=\max\limits_{{\mathbb R}^n}\rho_0(x)\omega_n
C^{(2+n)},$
 and $\omega_n$ is the volume of
a unit ball in ${\mathbb R}^n.$

Then the life span of the solution is finite.
\end{Theorem}

{\sc Remark 3.1} In \cite{Sideris} the system under consideration
is hyperbolic in physical space, therefore $n=3,$  $\alpha=1,$ and
$|R(t)|\le R(0)+\sigma t,$ the constant $\sigma=\left(
\frac{\partial P}{\partial \rho}|_{(\bar \rho,\bar S)}
\right)^{1/2}$ is the sound speed.

\medskip

Proof of Theorem 3.1. From the general Stokes formula we have
$m(t)=m(0),\,\eta(t)=\eta(0).$ The Jensen inequality together with
$\eta(0)\ge 0$ give us, according to \cite{Sideris},
$$\displaystyle\int\limits_{B(t)}P dx\ge (vol B(t))^{1-\gamma}
\left(\int\limits_{B(t)}\rho \exp \left(\frac{S}{\gamma}\right)
dx\right)^\gamma= $$
$$=(vol
B(t))^{1-\gamma} \left(\eta(0)+vol
B(t)\bar\rho\exp\left(\frac{\bar
S}{\gamma}\right)\right)^\gamma\ge \bar P vol B(t)=
\int\limits_{B(t)}\bar P dx.\eqno(32)$$

Let us note that if we integrate instead of ${\mathbb R}^n$ over
$B(t)$, we get an analog of Lemma 2.1 and Corollary 2.1. The only
difference will be in the  integral $I_{4}(t).$ Here $I_{4}(t)=-
\int\limits_{ S(t)} P ({\bf x,N}) dS, $ where $S(t)$ is the
boundary of $B(t)$. Therefore
$$\tilde G''(t)=F'(t)=
I_{1}(t)+I_{2}(t)+I_{3}(t)+I_{4}(t)\ge
I_{1,\phi}(t)+I_{2,\phi}(t)+n\int\limits_{B(t)}(P-\bar P) dx.$$

In this case $(\tilde G'(t))^2=F^2(t)\ge 2\tilde G(t)I_{1}(t) $
and $|M|^2\ge 2\tilde G(t)I_{2}(t).$ Thus, taking into account
(32), we have
$$F'(t)\ge \frac{F^2(t)+|{\bf M}|^2}{2\tilde G(t)}.\eqno(33)$$
 Further, since
$$\tilde G(t)\le\frac{1}{2}(C(1+t)^\alpha)^2(m(0)+\int\limits_{B(t)}\bar\rho
dx)= \frac{1}{2}(C(1+t)^\alpha)^2\int\limits_{B(t)}\rho_0 dx
=$$
$$\frac{1}{2}\max\limits_{{\mathbb R}^n}\rho_0(x)\omega_nC^{(2+n)}
(1+t)^{\alpha (2+n)},$$ we get from (33) that
$$F'(t)\ge
(\max\limits_{{\mathbb R}^n}\rho_0(x)\omega_nC^{(2+n)}
(1+t)^{\alpha (2+n)})^{-1} (F^2(t)+|{\bf M}|^2).$$

Integrating (34) we get the following.  For $|{\bf M}|=0$
$$-\frac{1}{F(t)}+
\frac{1}{F(0)}\ge
\frac{1}{A}\frac{(1+t)^{1-\alpha(2+n)}-1}{1-\alpha(2+n)},\quad
\alpha(2+n)\ne 1,\eqno(35)$$
$$-\frac{1}{F(t)}+
\frac{1}{F(0)}\ge \frac{1}{A}\ln (1+t),\quad \alpha(2+n)=
1.\eqno(36)$$ If $\alpha<\frac{1}{2+n},$ then we have from (35)
$$F(t)\ge \frac{AF(0)(1-\alpha(2+n))}{A(1-\alpha(2+n))-
F(0)((1+t)^{1-\alpha(2+n)}-1)}.$$ Respectively, from (36)we have
$$F(t)\ge \frac{AF(0)}{A- F(0)\ln(1+t)}.$$
Thus, if condition  (a) of Theorem 3.1 holds, then $F(t)$ become
infinite within a finite time, whereas it follows from the
Cauchy-Schwartz inequality that $F^2(t)\le 4 \tilde G (t)E_k(t),$
that is  $F(t)$ is finite at finite $t.$ We obtain a
contradiction.

Further, for $\alpha>\frac{1}{2+n}$  from (35) we get
that$$F(t)\ge \frac{AF(0)(\alpha(2+n)-1)(1+t)^{\alpha(2+n)-1}}
{F(0)-(F(0)-A(\alpha(2+n)-1))(1+t)^{\alpha(2+n)-1}},$$ it follows
that provided condition (c) of Theorem 3.1 holds, $F(t)$ goes to
infinity within a finite time.

Let us $|{\bf M}|\ne 0.$ Then
$$\arctan \frac{F(t)}{|{\bf M}|}\ge \arctan \frac{F(0)}{|{\bf
M}|}+\frac{|{\bf
M}|}{A}\frac{(1+t)^{1-\alpha(2+n)}-1}{1-\alpha(2+n)},\quad
\alpha(2+n)\ne 1,\eqno(37)$$
$$\arctan \frac{F(t)}{|{\bf M}|}\ge \arctan \frac{F(0)}{|{\bf
M}|}+\frac{|{\bf M}|}{A}\ln (1+t),\quad \alpha(2+n)= 1.\eqno(38)$$
One can conclude from (37) and (38), that if $\alpha\le
\frac{1}{2+n},$ then at any value of $F(0)$ in a finite time the
right-hand sides of these inequalities exceed $\frac{\pi}{2},$
whereas the left-hand sides are later then this number. This
contradiction shows that the solution cannot keep smoothness
provided condition (b) of Theorem holds.

At last, if $\alpha> \frac{1}{2+n},$ then an analogous
contradiction we can get from (37), if
$$\arctan
\frac{F(0)}{|{\bf M}|}+\frac{|{\bf
M}|}{A(\alpha(2+n)-1)}>\frac{\pi}{2},\eqno(39)$$ this results
conditions (d) and (e).

Thus, Theorem 3.1 is proved. $\Box$

\medskip

\medskip

{\sc Remark 3.1.} The statement of Theorem 1 from \cite{Sideris}
is a particular case of condition (c) of Theorem 3.1 for ${\bf
F}=0,\,\alpha=1, n=3.$ Let us analyze the condition (c) in the
physical space.
 For it we decompose the velocity  ${\bf V}$ into a sum
 of its radial and tangential components, ${\bf V}_r$ and ${\bf V}_\tau.$
 We mean that the  tangential component is a projection of velocity
 into subspace, orthogonal to the radius-vector ${\bf x}$, that is
$({\bf V}_\tau,{\bf x})=0.$ As the estimate
$$|F(0)|\le \max\limits_{{\mathbb
R}^n}\rho_0(x)\max\limits_{{\mathbb R}^n}|{\bf
V}_r(0,x)|C^{n+1}\omega_n$$ is true, then from condition (c) we
get that
$$(\alpha(2+n)-1)C<\max\limits_{{\mathbb R}^n}|{\bf V}_r(0,x)|,
$$
that is
$$4C<\max\limits_{{\mathbb R}^3}|{\bf V}_r(0,x)|.
$$

 It signifies that the initial perturbation of radial component of
 velocity (or the initial divergency, according to Remark 2.4)
 is sufficiently large comparing with the velocity  of the support
 expanding, $C,$ that it with the sound speed at infinity.

However, if $|{\bf M}|\ne 0,$  that is, according to Remark 2.4,
there exists an initial vorticity, then, as condition (d) shows,
initially the radial component can be later, since
$$|{\bf M}|\cot \frac {|{\bf
M}|}{A(\alpha(2+n)-1)}<A(\alpha(2+n)-1),$$ è
$$|{\bf M}|\cot \frac {|{\bf
M}|}{A(\alpha(2+n)-1)}\to A(\alpha(2+n)-1),\qquad |{\bf M}|\to
0.$$

Thus, condition (c) follows from (d) at $|{\bf M}|\to 0.$

At the same time, as condition (e) shows, the singularity may
appear due to a large vorticity only. Indeed,
$$|{\bf M}|\le\max\limits_{{\mathbb R}^n}\rho_0(x)\max\limits_{{\mathbb
R}^n}|{\bf V}_\tau(0,x)|C^{n+1}\omega_n,$$ it shows together with
 (e) that
 $$\pi(\alpha(2+n)-1)C<\max\limits_{{\mathbb R}^n}|{\bf V}_\tau(0,x)|,
$$
that is
$$4\pi C<\max\limits_{{\mathbb R}^3}|{\bf V}_\tau(0,x)|.
$$
However, at any case, it results the increase of the radial
component of velocity.

\medskip
\medskip\medskip

\section{Influence of damping to the singularity formation}

We assume that in the right hand side of (1) there arise an
additional forcing, namely, the dry friction with the coefficient
$\mu(t,x)$.  Thus, instead of (1), we consider
$$ \rho (\partial_t
{\bf V}+({\bf V},{\bf \nabla})\,{\bf V})+{\bf\nabla} P = {\bf
F}(t, x, \rho, {\bf V}, S, D^{|\alpha|}{\bf V})-\mu \rho{\bf
V}.\eqno(40)$$

In this situation the total energy is not conserved, however, this
function is non- increasing:
$$E(t)\le E(0).$$

It is known that if the density is initially compactly supported
then the dry friction (arbitrary large) do not prevent the
singularity formation, only delays it \cite {RozMag}.

At the same time, for the initial perturbation of nontrivial
constant state concentrated in bounded domain, the damping
prevents the singularity formation \cite {Wang}.

One can prove that if the friction is small, then   singularities
may arise both for solutions from class $\mathfrak K$ and for
compact perturbation of nontrivial constant state. We concentrate
at the first case, since the proof in the second one is analogous
(see also \cite {Wang} in this context).

Let us come back to the notation of Section 2, that is we will
write below $G(t), I_k(t), \,k=1,...,4$ instead of
$G_2(t),I_{k,\phi}(t),$ respectively.

We suppose that $|\mu(t,x)|\le \mu_0,$ with some positive constant
$\mu_0.$

Let us come back to the notation of Section 2, that is for
$\phi(|{\bf x}|)=\frac{1}{2}|{\bf x}|^2$ we will write below
$\,G(t), I_k(t),\,$ $ \,k=1,...,4$ instead of
$G_\phi(t),I_{k,\phi}(t),$ respectively.

Let us note that here Lemmas 2.2 and 2.3 are true, however, in the
expression for $G''(t)=F'(t)$ from Lemma 2.1 we have to add in
this new situation one more term. Namely, taking into account
Corollary 2.1 we can write this expression in the form
$$F'(t)= I_{1}(t)+I_{2}(t)+n(\gamma-1)E_p(t)-
\int\limits_{{\mathbb R}^n}\mu(t,x) ({\bf V},x)\rho dx.\eqno(41)$$
If $\mu$ is constant, then $$\int\limits_{{\mathbb R}^n}\mu(t,x)
({\bf V},x)\rho dx=\mu F(t),$$ in general case, using the estimate
$$|\int\limits_{{\mathbb R}^n}\mu (t,x)({\bf
V},x)\rho dx|\le 2\mu_0\sqrt{E(0)G(t)},$$ we obtain only the
inequality
$$F'(t)\ge I_{1}(t)+I_{2}(t)+n(\gamma-1)E_p(t)-2\mu_0\sqrt{E(0)G(t)}.
\eqno(42)$$

The value of  $|{\bf M}(t)|$ is not conserved here. If $\mu$ is
constant, then
$$M_k(t)\le M_k(0) e^{-\mu t},\,k=1,...,K.$$
In general case we note that $(M'_k(t))^2\le 4\mu_0^2E_k(t)G(t)$
and apply estimate (11). Integrating the two-sided inequality for
$M'_k(t)$ we obtain that
$$M_k(0)-\mu_0((\sqrt{E(0)}t+\sqrt{G(0)})^2-G(0))\le
M_k(t)\le $$$$M_k(0)+
\mu_0((\sqrt{E(0)}t+\sqrt{G(0)})^2-G(0)).\eqno(43)$$

\medskip

\begin{Theorem}
For sufficiently small  $\mu_0$ there are initial data (6') such
that the solution to the Cauchy problem (40),(2),(5),(6') cannot
belong to the class $\mathfrak K$  for all $t\ge 0.$
\end{Theorem}

Proof. As follows from (41), (8), (10)
$$F'(t)\ge
\frac{F^2(t)+|{\bf
M}(t)|^2}{2G(t)}+\frac{n(\gamma-1)C}{(G(t))^{\frac{(\gamma-1)n}{2}}}
-2\mu_0 \sqrt{E(0)G(t)}.
$$

For $\gamma\le 1+\frac{2}{n}$ it results inequality
$$
F'(t)\ge \frac{
F^2(t)+L^2\Psi(\mu_0,t)}{2(\sqrt{E(0)}t+\sqrt{G(0)})^{2}}.\eqno(44)$$
Here we use the notation of Theorem 2.1 ( $|{\bf M}|$ denotes now
$|{\bf M}(0)|$ in the expression for
 $L$), and we introduce a function
$$\Psi(\mu_0,t)=1-\frac{1}{L^2}(4\mu_0\sqrt{E(0)}(\sqrt{E(0)}t+
\sqrt{G(0)})^3 +$$
$$n E(0)\mu_0^2t^2(\sqrt{E(0)}t+2\sqrt{G(0)})^2+2\mu_0\sqrt{n}|{\bf
M}(0)|\sqrt{E(0)}t(\sqrt{E(0)}t+2\sqrt{G(0)})).$$

 Let us fix a
positive constant $\Psi_*^2<1.$ Choosing $\mu_0$ sufficiently
small we can obtain that the inequality $\Psi(t)>\Psi^2_*$ will
hold for all $t\in [0,T(\mu_0)),\,$ moreover, $T(\mu_0)\to \infty
$ as $\mu_0\to 0.$ Thus, for such $t$ we have from (44) that
$$
F'(t)\ge \frac{
F^2(t)+L^2\Psi_*^2}{2(\sqrt{E(0)}t+\sqrt{G(0)})^{2}},$$ the
integration results
$$\arctan \frac{F(t)}{L\Psi_*}\ge
\arctan \frac{F(0)}{L\Psi_*}
+\frac{L\Psi_*}{2\sqrt{E(0)}}\left(\frac{1}{\sqrt{G(0)}}-
\frac{1}{\sqrt{E(0)}t+\sqrt{G(0)}}\right).$$

Let us denote $t_*$ a unique positive solution (in $t$) to
equation
$$\arctan
\frac{F(0)}{L\Psi_*}
+\frac{L\Psi_*}{2\sqrt{E(0)}}\left(\frac{1}{\sqrt{G(0)}}-
\frac{1}{\sqrt{E(0)}t+\sqrt{G(0)}}\right)=\frac{\pi}{2},$$ which
always exists if $$F(0)\ge {L\Psi_*}
\cot\frac{L\Psi_*}{2\sqrt{E(0)G(0)}}.$$ It suffices to choose
$\mu_0$ such small that $T(\mu_0)>t_*.$

 The
case $\gamma>1+\frac{2}{n}$ can be treated analogously.

Thus, Theorem 5.1 is proved. $\Box$

\medskip

{\sc Remark 4.1.} The function
${L\Psi_*}\cot\frac{L\Psi_*}{2\sqrt{E(0)G(0)}}$ is decreasing in
$\Psi_*$ for  $0<\Psi_*\le 1,$ therefore
$${L\Psi_*}\cot\frac{L\Psi_*}{2\sqrt{E(0)G(0)}}>
{L}\cot\frac{L}{2\sqrt{E(0)G(0)}},$$ and if condition (14) is
satisfied, the condition
$$F(0)\ge {L\Psi_*}\cot\frac{L\Psi_*}{2\sqrt{E(0)G(0)}},$$
sufficient for the formation of singularity in the system with
damping, generally speaking, does not hold. In this sense the
damping prevents the singularity formation.

\medskip

\section{Influence of rotation to the singularity formation}

In many meteorological problems for the modelling of the
atmospherical processes they use systems of equations, analogous
to (1--3). However, in this problems the air motion is considered
under rotating Earth, therefore it needs to take into account the
Coriolis force. It is defined as ${\bf V}\times 2{\bf \Omega},$
where ${\bf \Omega}$ is a constant vector of the Earth angular
rotation (see, for example, \cite{Pedlosky},\cite{Gill}).

The question on sufficient conditions of singularity formation
from smooth initial data is very important for the meteorology,
where the discontinuity is associated with the atmospherical
front.

The vertical scale of atmospherical processes, considered in the
problems of weather forecast, is small compared with the
horizontal one (no more then 10 km in vertical against several
hundred kilometers by horizontal), therefore equations describing
horizontal and vertical processes, are not equal in rights. In the
vertical direction the so called hydrostatic approximation is
usually assumed. Its general sense is that vertical velocity and
acceleration are assumed to be small with respect to horizontal
ones (see, for example, \cite{Pedlosky},\cite{Gill}). (This
approach is not acceptable, of course, for a description of
 small scales processes  such that spouts, convection near frontal zones,
typhoons generation, and so on.)

However, if the scale is such that the assumption on a smallness
   of vertical processes is acceptable, then
 it is convenient to average the
 primitive tree dimensional system of equation over hight
and to deal with a simpler two dimensional in space system of
 equation where both space equations are already equal
 in rights
\cite{Obukhov},\cite{Alishaev}. Vertical processes are hidden now,
however, they are not eliminated  from consideration. One can see
it, for example, in the fact that the value of adiabatic exponent
changes in the averaged over high system.

 We also restrict ourself to the case $n=2.$
\medskip

Thus, instead of (1) we consider
$$ \rho (\partial_t
{\bf V}+({\bf V},{\bf \nabla})\,{\bf V})+{\bf\nabla} P = {\bf
F}(t, x, \rho, {\bf V}, S, D^{|\alpha|}{\bf V})+l\rho{\bf
V}_\bot.\eqno(45)$$ We  denote ${\bf V}_\bot(v_\bot^1,v_\bot^2),$
a vector with components $\,v_{\bot}^i=e_{.j}^i v^j,\,i,j=1,2$
where $ e_{ij} $ is the Levi-Civita tensor,  $l =l(x)$ is the
Coriolis parameter.

\medskip

We consider as before the solutions from the class $\mathfrak K.$

Among compactly supported solutions for $l=0$ there exist no
globally smooth ones \cite{Makino}.

In contrast to the case for $l\ne 0$ we can construct explicitly a
stationary nontrivial compactly supported solution.

Let us  ${\bf F=0}$ and $l$ is constant. We consider the
isentropic gas, that is the state equation is $P=A
\rho^{\gamma},\,A=const)$. We will seek a solution of the form
$(\rho,{\bf V}, P),$ where
$$ {\bf V}= f(\theta){\bf r}_{\bot},\,
\rho=(g(\theta))^{1/(\gamma-1)}, \, P=A \rho^{\gamma}.$$ Here
$\theta=|{\bf x}|^2/2\, ,\,$ $f(\theta) $ is an arbitrary smooth
function supported on a segment $[a,b]\in [0,\infty),\,$ $$
g(\theta) = C+\frac{1}{2 K}\int_0^\theta (f^2(\xi)-l f(\xi))
d\xi,\,$$ with the constant $K=\frac{A\gamma}{\gamma-1}$. We can
always choose the constant $C$ such that $ g(\theta) $ will vanish
as $|{\bf x}|\to \infty.$ For example, if $
f(\theta)=\frac{l}{\nu}(\nu-\theta)$ for $0\le\theta\le\nu<\infty
$ and $f(\theta)=0$ for $\theta>\nu$ we have $
g(\theta)=\frac{l^2}{12 K \nu^2}(\nu^3+\theta^2(2\theta-3\nu)),\,$
$ 0\le\theta\le\nu<\infty,$ and $g(\theta)=0,\,\theta>\nu.$

These solutions correspond to compactly supported stationary
divergence-free flows.

For the constant Coriolis parameter we can also construct a
non-stationary periodic in time global solution from the class
$\mathfrak K$ with linear profile of velocity $${\bf
V}=\alpha(t){\bf x}+\beta(t){\bf x}_\bot,$$ acting in the spirit
of Subsection 2.1 (see \cite{RozPas}, \cite{RozJMS},
\cite{RozNova}for detail).

\medskip

However, in the rotational case, too, there exists initial data,
resulting in a singularity formation in the class of solutions
with a finite moment of mass and total energy.

If $|l(x)|\le l_0,$ these data can be found exactly as in the case
of small damping, described in Section 4. It is sufficient to note
that
$$|\int\limits_{{\mathbb R}^2}l(x) ({\bf V}_\bot,x)\rho dx|\le
2l_0\sqrt{E(0)G(t)},$$ and to proceed as before changing $\mu$ to
$l.$ The relative conclusion is that if the rotation is small we
can still find  initial conditions resulting singularity
satisfying to an analog of (14).

 It will be convenient for us to consider $l$ as a
constant.

\medskip

{\sc Remark 5.1} In \cite{LiuTadmor} the question  was
investigated whether the rotation prevents the singularity
formation for the pressureless gas dynamics. Here the basic
equation is the nonlinear transport equation with rotational
forcing, namely
$$ \partial_t
{\bf V}+({\bf V},{\bf \nabla})\,{\bf V}= l\rho{\bf V}_\bot.$$ The
answer is "conditionally yes". Specifically, in the case $l\ne 0$
the solution is globally smooth if and only if
$$
\forall x\in {\mathbb R}^2 \quad 2l\omega_0(x)+\eta_0(x)^2<l^2,
$$
where $\omega_0(x)=(V_{02}(x))'_{x_1}-(V_{01}(x))'_{x_2},\quad
\eta_0(x)=\lambda_2(0)-\lambda_1(0),\quad \lambda_j(0)$ are
eigenvalues of Jacobian of initial velocity.

Let us point out that this condition is violated if in a certain
point $x_0$ the value of $\omega_0(x_0)>>1$ (meteorologists would
say that near $x_0$ there exists a significant cyclonic vorticity;
in the northern hemisphere, where the Coriolis parameter  $l$ is
positive, the cyclonic rotation is anticlockwise).

Let us note also that we cannot do here the limit pass as $l\to
0,$ since if $l=0,$ then the solution is globally smooth if and
only if $\lambda_j(0)\ge 0,\,j=1,2.$

\medskip

It is especially interesting that the presence of rotation in the
real gas dynamic is in some sense convenient for the singularity
formation. More precisely, there are situations where only
significant initial vorticity provokes the singularity at any
initial divergency. As we have seen, it is impossible in the
rotation free case, where the divergency must be significantly
positive.

\medskip

In in the rotational case the mass $m$ is conserved, the total
energy $E(t)$ is non-increasing (it is conserved for ${\bf F=0}$).
Corollary 2.1 and inequalities (7--11), estimates (12) and (13)
are true. The angular momentum balance law has now the form
$${\mathcal M}=lG(t)+F_\bot(t)=const,\eqno(46)$$
where $$F_\bot(t)=\int\limits_{{\mathbb R}^2}({\bf V}_\bot,x)\rho
dx.$$

\medskip

In this new situation the following Lemma holds:

\begin{Lemma} For the solution of class $\mathfrak K$ to
system (45), (2), (5)  equalities
$$G'(t)=F(t),\eqno(47)$$
$$F'(t)=I_1(t)+I_2(t)+I_3(t)+lF_\bot(t),\eqno(48)$$
$$F_\bot'(t)=-l F(t),\eqno(49)$$
$$G''(t)+l^2G(t)=2(2-\gamma)E_k(t)+\Theta_\gamma(t),\eqno(50)$$
where $\Theta_\gamma(t)=2(\gamma-1)E(t)+l{\mathcal M},$ take
place.
\end{Lemma}

\medskip

The proof of Lemma 5.1 is absolutely analogous to the Lemma 2.1
proof, the general Stokes formula in this two-dimensional case
looks like the Green' formula. $\Box$

\medskip

{\sc  Remark 5.2.} The balance law (46) follows from (47) and
(49).

\medskip

A new important circumstance in this situation is a boundedness of
$G(t)$ above.

\begin{Lemma}
For the solution of class $\mathfrak K$ to system (45), (2), (5)
the inequality
$$0< G_-\le G(t) \le G_+\eqno(51)$$
takes place, where $G_-$ and $G_+$ are positive constants. For
example, on can take
$$G_-=\left(\frac{C}{E(0)}\right)^{1/(\gamma-1)},\quad
G_+=\frac{1}{l^2}(\sqrt{\Theta_2(0)} +\sqrt{2E(0)})^2.$$
\end{Lemma}

\medskip

{\sc Remark 5.3.} Let us denote $J_i=\int\limits_{{\mathbb
R}^2}V_i\rho dx,\,i=1,2.$ The  value $J_1^2(t)+J_2^2(t)$ is
conserved for ${\bf F=0}$ (for example, \cite{Gordin}). The
estimates (51) can be refined (see \cite{RozFAO}), if this
conserved values is positive.

\medskip

Proof of Lemma 6.2. It follows from (50) that
$$G''(t)+l^2G(t)\le\Theta_2(0),$$
however due to the resonance phenomenon we cannot prove the
boundedness of the solution to this inequality without taking into
account additional properties of $G(t).$

Suppose that there exists a point $t_1$ such that $ G'(t_1)>0$ (
otherwise, $G(0)$ is the upper bound of $G(t)$). Let us denote
$\epsilon(t_1)=G(t_1)-\frac{\Theta_2(0)}{l^2}.$ If for all $t_1$
it occurs that $\,\epsilon(t_1)\le 0,$  then the upper bound of
$G(t)$ is $\frac{\Theta_2(0)}{l^2}.$ Let us suppose that one can
find $t_1$ such that $\,\epsilon:=\epsilon(t_1)> 0.$ Then there
exists such $\tau>0,$ that $G(t)\ge
\frac{\Theta_2(0)}{l^2}+\epsilon$ at $t\in [t_1, t_1+\tau).$

Thus, it follows from (52) that at   $t\in [t_1, t_1+\tau)$
$$G''(t)\le -\epsilon l^2.$$
We integrate this inequality twice with respect to $t$ from $t_1$
to $t\in [t_1, t_1+\tau),$ and we get that
$$G(t)\le -\epsilon\frac{l^2}{2}(t-t_1)^2+G'(t_1)(t-t_1)+G(t_1):=G^+(t).$$
The quadratic function $G^+(t)$ is maximal at
$t=t_*=t_1+\frac{G'(t_1)}{\epsilon l^2}>t_1.$  Therefore $G(t)\le
G^+(t_*) $ at $t_*\le t_1+\tau,$ and $G(t)\le
G^+(t_1+\tau)<G^+(t_*), $ at  $t_*> t_1+\tau.$ At any case, taking
into account  (7), we get that
$$G(t)\le G^+(t_*)= G(t_1)+ \frac{(G'(t_1))^2}{2\epsilon
l^2}\le G(t_1) + \frac{4EG(t_1)}{2\epsilon
l^2}=$$$$=\frac{\Theta_2(0)}{l^2}+\epsilon+\frac{4E(\frac{\Theta_2(0)}
{l^2}+\epsilon)}{2\epsilon
l^2}=\frac{\Theta_2(0)+2E}{l^2}+\epsilon+
\frac{2E\Theta_2(0)}{l^4\epsilon}.$$ Minimizing the right hand
side in $\epsilon$ we get the estimation of $G(t)$ from above.

The lower bound indicated in the Lemma statement can be obtained
from Lemma 2.3. $\Box$

\begin{Corollary}
In the situation of Lemma 6.2 $|F(t)|$ is bounded for all $t>0.$
\end{Corollary}

Indeed, it follows from (7) and Lemma 6.2. $\Box$

\medskip

Let us denote ${\mathcal K}={\mathcal M}^2-l^2G_+^2+\delta ,$
where $\delta=CG_-^{1-\frac{(\gamma-1)n}{2}},$ if $\gamma\le
1+\frac{2}{n},$ and $\delta=CG_+^{1-\frac{(\gamma-1)n}{2}},$
otherwise. The following Theorem takes place.

\begin{Theorem}
Solutions to (45), (2), (5) cannot belong to the  class $\mathfrak
K$ for all $t\ge 0$ if the initial data are such that
$${\mathcal K}>0,\eqno(53)$$
or
$${\mathcal K}\le 0,\quad F(0)>\sqrt{-{\mathcal K}}\eqno(54)$$
hold.
\end{Theorem}

Proof of the Theorem. It follows from (8), (10), (13), (46), (48)
and Lemma 2.3 that
$$F'(t)\ge
\frac{F^2(t)+F_\bot^2(t)}{2G(t)}+\frac{C}{G^{(\gamma-1)n/2}}+
l{\mathcal M}-l^2G(t)\ge$$$$\ge \frac{F^2(t)+{\mathcal
K}}{2G_+}.\eqno(55)$$ We integrate (55) and see that if (53) or
(54) are satisfied, then $F(t)$ become unbounded in a finite time.
This contradicts to Corollary 6.1. $\Box$

\medskip
\medskip

{\sc Remark 5.4.} The estimate from above of the time of the
singularity formation can be easily obtained from inequality (55).

\medskip

{\sc Remark 5.5.} Let us analyze, for example, condition (53),
that is
$$({\mathcal M}-lG_+)({\mathcal M}+lG_+)>-\delta.\eqno(56)$$
Taking into account the expression for $G_+,$ given by Lemma 5.2,
we can  re-write inequality (56) as
$$(\Theta_2(0)-l{\mathcal M}+\sqrt{2E\Theta_2(0)})(\Theta_2(0)+2\sqrt{2E
\Theta_2(0)})<\frac{\delta}{4}l^2.$$ Let us note that since
$\Theta_2(0)>0,$ (56) does not hold for $\delta=0,$ in the case of
pressure free gas dynamics.

If $\delta>0,$ inequality (56) is satisfied for
$\Theta_2(0)<\Theta_*(\delta,l),$ with a constant
$\Theta_*(\delta,l).$

Thus, coming back to inequality (53), we see that it is true, if
$l{\mathcal M}<\Theta_*(\delta,l)-2E,$ that is
$$lF_\bot(0)<l^2G(0)-2E+\Theta_*(\delta,l).$$
For large $E$ and small (however not equal to zero!) $|l|$ (this
is in the real meteorological situation) the value of $F_\bot(0)$
is negative. Since
$$F_\bot(t)=-\frac{1}{2}\int\limits_{{\mathbb R}^2}((\rho V_2)'_{x_1}-
(\rho V_1)'_ {x_2}) dx,$$
 it signifies that there exists initially a cyclonic vorticity.
 Possibly, this observation would help to explain
 the well known for meteorologists fact that in the extratropical
 zone always inside of cyclone  an atmospherical front exists.

\medskip

{\sc Remark 5.6} The results of this section can be with
respective modification extended to the case of non-constant
Coriolis parameter, which, however, differs little from a
constant.

\medskip

{\sc Remark 5.7}  One can consider in the rotational case, too,
the viscid term of form (30). Then we can obtain as a corollary
from Theorem 5.1 the following result. A solution from the class
$\mathfrak K$ with the conditions of decay at infinity for the
velocity and its derivatives (see Remark 2.8), loses the initial
smoothness provided  (53) or (54) hold.

\medskip

{\sc Remark 5.8}  A certain result concerning sufficient condition
for the singularity formation for the Euler equations on a
rotating plane demonstrating another approach can be found
in\cite{RozDU}. I is possible to consider other exterior forces
then ones mentioned in this paper. For example, \cite{RozZur}
deals with sufficient conditions of the smoothness loss for
solutions to the gas dynamic equation with exterior force of
geopotential type (besides the Coriolis force).

\end{document}